\documentclass[12pt]{article}
\usepackage{mathtools} 
\usepackage{amsmath, amssymb} 
\usepackage{amsthm}
\usepackage{enumerate}  
\usepackage{url} 
\usepackage{authblk}
\usepackage{tikz}

\usepackage{graphicx}
\usepackage{caption,subcaption}

\newcommand\cx{{\mathbb C}}

\newcommand\ints{{\mathbb Z}}
\newcommand\re{{\mathbb R}}
\newcommand\rats{{\mathbb Q}}

\DeclarePairedDelimiter\abs{\lvert}{\rvert}%
\DeclarePairedDelimiter\norm{\lVert}{\rVert}%

\makeatletter
\let\oldabs\abs
\def\abs{\@ifstar{\oldabs}{\oldabs*}}
\let\oldnorm\norm
\def\norm{\@ifstar{\oldnorm}{\oldnorm*}}
\makeatother

\newcommand\opk[1]{\mathop{\mathrm{#1}}\nolimits}

\newcommand\comp[1]{{\mkern2mu\overline{\mkern-2mu#1}}}
\newcommand\pmat[1]{\begin{pmatrix} #1 \end{pmatrix}}
\newcommand\seq[4]{#1_{#2},#1_{#3},\ldots,#1_{#4}}

\newtheoremstyle{plainsl}%
	{\topsep}
	{\topsep}
	{\slshape} 
	{}
	{\normalfont\bfseries}
	{.}
	{ }
	{}

\swapnumbers

{\theoremstyle{plainsl}
\newtheorem{theorem}{Theorem}[section]
\newtheorem{lemma}[theorem]{Lemma}
\newtheorem{corollary}[theorem]{Corollary}}
{\theoremstyle{remark}
}

\renewcommand\proof{\noindent\textsl{Proof. }}
\newcommand\sqr[2]{{\vbox{\hrule height.#2pt
    \hbox{\vrule width.#2pt height#1pt \kern#1pt
        \vrule width.#2pt}\hrule height.#2pt}}}
\renewcommand\qed{%
	\ifmmode\eqno\sqr53
	\else\nolinebreak\ \hfill\sqr53\medbreak\fi}

\DeclareMathOperator{\col}{col}

\newcommand\ip[2]{\langle#1,#2\rangle}
\newcommand\one{{\bf1}}

\newcommand\sym[1]{\opk{Sym}(#1)}

\usepackage{blkarray}

\title{Pretty good state transfer in discrete-time quantum walks}
\date{}

\author{Ada Chan}

\author{Hanmeng Zhan}

\affil{Department of Mathematics and Statistics, York University, Toronto, ON, Canada\\\texttt{\{ssachan, h3zhan\}@yorku.ca}}

\begin{document}
\maketitle

\begin{abstract}
We establish the theory for pretty good state transfer in discrete-time quantum walks. For a class of walks, we show that pretty good state transfer is characterized by the spectrum of certain Hermitian adjacency matrix of the graph; more specifically, the vertices involved in pretty good state transfer must be $m$-strongly cospectral relative to this matrix, and the arccosines of its eigenvalues must satisfy some number theoretic conditions. Using normalized adjacency matrices, cyclic covers, and the theory on linear relations between geodetic angles, we construct several infinite families of walks that exhibits this phenomenon. 
\end{abstract}

\section{Introduction}

A quantum walk with nice transport properties is desirable in quantum computation. For example, Grover's search algorithm \cite{Grover1996} is a quantum walk on the looped complete graph, which sends the all-ones vector to a vector that almost ``concentrate on" a vertex. In this paper, we study a slightly stronger notion of state transfer, where the target state can be approximated with arbitrary precision. This is called pretty good state transfer.

Pretty good state transfer has been extensively studied in continuous-time quantum walks, especially on the paths \cite{GodsilKirkland2012, Vinet2012,Banchi2017,Coutinho2017,VanBommel2019}. However, not much is known for the discrete-time analogues. The biggest difference between these two models is that in a continuous-time quantum walk, the evolution is completely determined by the adjacency or Laplacian matrix of the graph, while in a discrete-time quantum walk, the transition matrix depends on more than just the graph - usually, it is a product of two unitary matrices:
\[U = S C,\]
where $S$ permutes the arcs of the graph, and $C$, called the \textsl{coin matrix}, send each arc to a linear combination of the outgoing arcs of the same vertex. Thus, the relation between the walk and the graph can be rather complicated.

In this paper, we consider a discrete-time quantum walk where the spectrum of $U$ is determined by certain Hermitian adjacency matrix $H$ of the graph. A generalization of this walk, called the twisted Szegedy walk, was proposed in \cite{Higuchi2014}, and applied in \cite{Segawa2020,Kubota2021}. We show that pretty good state transfer in this model reduces to $m$-strong cospectrality relative to $H$ and some number theoretic conditions on the arccosines of its eigenvalues. Using normalized adjacency matrices, cyclic covers, and the theory on linear relations between geodetic angles, we construct several infinite families of walks that exhibits this phenomenon.

\section{A discrete-time quantum walk \label{sec:dqw}}

Let $X$ be a connected undirected graph. Let $W$ be a complex adjacency matrix of $X$, that is, $W=(w_{ab})_{ab}$ where $w_{ab}\ne 0$ if and only if $a$ is adjacent to $b$. We say $W$ is \textsl{normalized} if the entrywise product $W\circ \comp{W}$ is row-stochastic, that is, 
\[(W \circ \comp{W} )\one = \one.\]

We construct a quantum walk with respect to a given pair $(X, W)$ of graph and normalized complex adjacency matrix. To each edge $\{a, b\}$ of $X$, we associate a pair of opposite arcs, $(a, b)$ and $(b,a)$; we call $a$ and $b$ the \textsl{tail} and the \textsl{head} of the arc $(a, b)$, respectively.  A \textsl{quantum state} is a complex-valued function on the arcs of $X$. Over time, it evolves according to a unitary \textsl{transition matrix} $U$, which we define now. Let $R$ be the permutation matrix that sends each arc $(a,b)$ to its opposite direction $(b,a)$. Let $N_t$ be the weighted tail-arc incidence matrix, where
\[(N_t)_{u, (a, b)} = \begin{cases}
w_{ab},\quad u=a,\\
0,\quad u\ne a.
\end{cases}\]
Since $W$ is normalized, $N_tN_t^*=I$ and $N_t^*N_t$ is a projection. The transition matrix $U$ is then given by
\[U := R (2N_t^* N_t - I).\]
After $t$ steps, the quantum walk will be in the state $U^t x$, if it started with the state $x$.

\section{Spectral decomposition of $U$ \label{sec:spd}}
The spectral decomposition of a product of two involutions has been studied in \cite{Szegedy2004}, and widely applied to quantum walks whose transition matrices satisfy this property \cite{Kendon2003, Szegedy2004,Higuchi2014, Zhan2020, Kubota2021}. In this section, we consider the quantum walk with respect to a pair $(X, W)$ of graph and normalized complex adjacency matrix, and let $U$ be the associated transition matrix:
\[U = R (2N_t^* N_t - I).\]
We establish a correspondence between the spectrum of $U$ and that of the Hermitian adjacency matrix
\[H:= W \circ W^*.\]

\begin{lemma}\label{lem:H}
We have
\[N_t^* R N_t = H.\]
\end{lemma}
\proof
 Pick any two vertices $a$ and $b$. The $ab$-entry of $N_t R N_t^*$ is 
 \begin{align*}
 e_a^T N_t R N_t^* e_b
 &= \left(\sum_{u \sim a} w_{au} e_{(a,u)}^T\right) R \left(\sum_{v \sim b} \comp{w_{bv}} e_{(b,v)}\right)\\
 &= \sum_{u\sim a} \sum_{v \sim b} w_{au} \comp{w_{bv}}R_{(a,u), (b,v)}\\
 &=\begin{cases}
 w_{ab} \comp{w_{ba}},\quad a\sim b,\\
 0,\quad a\not\sim b,
 \end{cases}
 \end{align*}
 where the last equality follows as $R$ is the permutation matrix that reverses each arc.
 \qed

Note that 
\[R^2=(2N_t^*N_t-I)^2=I,\] 
and so $U$ is indeed a product of two involutions. As every involution equals twice a projection minus the identity, the following result turns out to be useful in determining the eigenspaces of $U$.

\begin{lemma}\cite[Ch 2]{Zhan2018} \label{lem:spd}
Let $K$ and $L$ be matrices with orthonormal columns. Let $P=KK^*$, $Q=LL^*$ and $S=L^*K$. Let
\[ U = (2P-I)(2Q-I).\]
Then the eigenvalues of $SS^*$ lie in $[0,1]$. Moreover, we have the following spectral correspondence between $S$ and $U$.
\begin{enumerate}[(i)]
\item The $1$-eigenspace of $U$ is the direct sum
\[(\col(P)\cap \col(Q)) \oplus (\ker(P)\cap\ker(Q)).\]
Moreover, the map $z\mapsto Lz$ is an isomorphism from the $1$-eigenspace of $SS^*$ to $\col(P)\cap \col(Q)$.
\item The $(-1)$-eigenspace of $U$ is the direct sum
\[(\col(P)\cap \ker(Q)) \oplus (\ker(P)\cap\col(Q)).\]
Moreover, $z\mapsto Kz$ is an isomorphism from $\ker(S)$ to $\col(P)\cap\ker(Q)$, and $z\mapsto Lz$ is an isomorphism from $\ker(S^*)$ to $\ker(P)\cap\col(Q)$.
\item The remaining eigenspaces of $U$ are completely determined by the eigenspaces of $SS^*$. To be more specific, let $\mu\in(0,1)$ be an eigenvalue of $SS^*$ and let $\theta \in \re$ be such that $\cos(\theta)=2\mu-1$. The map
\[z\mapsto ((\cos(\theta)+1)I - (e^{i\theta}+1)P)Lz\]
is an isomorphism from the $\mu$-eigenspace of $SS^*$ to the $e^{i\theta}$-eigenspace of $U$, and the map
\[z\mapsto ((\cos(\theta)+1)I - (e^{-i\theta}+1)P)Lz\]
is an isomorphism from the $\mu$-eigenspace of $SS^*$ to the $e^{-i\theta}$-eigenspace of $U$.
\end{enumerate}
\end{lemma}

We now apply the above lemma to our walk, with the observation that $L=N_t^*$ and $2SS^*-I =N_t^*RN_t$.

\begin{theorem}\label{thm:eprojs}
Let $U$ be the transition matrix with respect to $(X,W)$, and let $H=W\circ W^*$. For each eigenvalue $e^{i\theta}$ of $U$, there is an eigenvalue $\lambda$ of $H$ such that $\lambda=\cos(\theta)$. Moreover, if $F_{\theta}$ denotes the $e^{i\theta}$-eigenprojection, and $E_{\lambda}$ denotes the $\lambda$-eigenprojection, then $F_{\theta}$ and $E_{\lambda}$ are related in the following way.
\begin{enumerate}[(i)]
\item If $\theta=0$, then $\lambda=1$ and
\[ N_t F_0 N_t^* = E_1.\]
\item If $\theta=\pi$, then $\lambda=-1$ and
\[N_t F_{\pi}N_t^* = E_{-1}.\]
\item If $\theta \in(-\pi,0)\cup(0,\pi)$, then
\[N_t F_{\theta}N_t^* = \frac{1}{2}E_{\lambda}.\]
\end{enumerate}
\proof
Note that with $L=N_t^*$ in Lemma \ref{lem:spd}, we have
\[\ker(Q) = \ker(N_t^*N_t)=\ker(N_t).\]
Suppose $\theta=0$. By Lemma \ref{lem:spd} (i), the projection onto $\col(P)\cap \col(Q)$ is given by 
\[F_0= N_t^* E_1 N_t,\]
and so $N_t F_0 N_t^*=E_1$. Now suppose $\theta=\pi$. By Lemma \ref{lem:spd}, the map $z\mapsto N_t^* z$ defines an isomorphism from $\ker(S^*)$ to $\ker(P)\cap\col(Q)$. Since $\ker(S^*) = \ker(SS^*)$ and $2SS^*-I=H$, we see that
\[F_1 = N_t^*E_{-1}N_t.\]
Finally, suppose $\theta\in(-\pi,0)\cup(0,\pi)$. Lemma \ref{lem:spd} (iii) and a bit of calculation show that
\[F_{\theta} = \frac{1}{2\sin^2(\theta)}(N_t^*- e^{i\theta}R N_t^*)E_{\lambda}(N_t - e^{-i\theta} N_t R).\]
Therefore,
\[N_tF_{\theta}N_t^* = \frac{1}{2\sin^2(\theta)}(I-e^{i\theta}H) E_{\lambda}(I-e^{-i\theta}H) = \frac{1}{2}E_{\lambda}.\tag*{\sqr53}\]
\end{theorem}

\begin{figure}[h]
\centering
\begin{minipage}{0.4\textwidth}
\includegraphics[width=4cm]{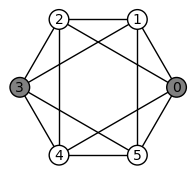}
\end{minipage}%
\begin{minipage}{0.6\textwidth}
\begin{itemize}
\item $W=\frac{1}{2}A(X)$
\item $N_t^*e_a=\frac{1}{2}(e_{(0,1)} + e_{(0,2)} + e_{(0,4)}+e_{(0,5)})$
\item $N_t^* e_b=\frac{1}{2}(e_{(3,1)} + e_{(3,2)} + e_{(3,4)}+e_{(3,5)})$
\item $t=6$
\end{itemize}
\end{minipage}
\caption{A graph with perfect state transfer}
\label{fig:pst}
\end{figure}

We are interested in transferring a state that ``concentrate" on some vertex $a$, that is, a state whose support is a subset of the outgoing arcs of $a$. One natural candidate for such states is $N_t^* e_a$, where $e_a$ is a the characteristic vector of the vertex $a$. In \cite{Zhan2019}, it is shown that for certain family of graphs, there are vertices $a$ and $b$, an integer time $t$ and a complex number $\gamma$ such that
\[U^t (N_t^* e_a) = \gamma (N_t^* e_b).\]
This is called \textsl{perfect state transfer} at time $t$. Figure \ref{fig:pst} illustrates one such example.

Unfortunately, perfect state transfer is rare. The main goal of this paper is to characterize a relaxation called \textsl{pretty good state transfer}, which occurs if there is a unimodular complex number $\gamma$ such that for any $\epsilon>0$, there is an integer time $t$ with
\[\norm{U^t N_t^*e_a-\gamma N_t^*e_b}<\epsilon.\]
The following result comes in handy when we analyze both perfect state transfer and pretty good state transfer from $N_t^* e_a$ to $N_t^*e_b$.

\begin{lemma}\label{lem:paral}
Let $e^{i\theta}$ be an eigenvalue of $U$ with eigenprojection $F_{\theta}$. Let $\lambda=\cos(\theta)$. Let $E_{\lambda}$ be projection onto the $\lambda$-eigenspace of $H$. Then 
\[F_{\theta} N_t^* e_a = \alpha F_{\theta}N_t^* e_b\]
if and only if 
\[E_{\lambda} e_a = \alpha E_{\lambda} e_b.\]
\end{lemma}
\proof
The identity
\[F_{\theta} N_t^* e_a = \alpha F_{\theta}N_t^* e_b\]
is equivalent to
\begin{equation}
F_{\theta} N_t^*( e_a- \alpha e_b)=0.\label{eqn:paral}
\end{equation}
Since $F_{\theta}$ is a projection, $N_t F_{\theta}N_t^*$ is positive semidefinite. Thus, \eqref{eqn:paral} holds if and only if 
\[N_t F_{\theta} N_t^*( e_a- \alpha e_b)=0.\]
The result now follows from Theorem \ref{thm:eprojs}.
\qed

\section{Perfect and pretty good state transfer in a generic quantum walk}
We develop the theory for perfect and pretty good state transfer in a generic quantum walk. Let $U$ be a unitary matrix acting on $\cx^n$. Consider two states $x$ and $y$ in $\cx^n$. We say there is \textsl{perfect state transfer} from $x$ to $y$, if there exists $t\in \ints$ and $\gamma \in \cx$ such that
\[U^t x = \gamma y.\]
Clearly, the above equation implies $\abs{\ip{U^t x}{y}}=1$; on the other hand, by Cauchy-Schwarz inequality, 
\[\abs{\ip{U^t x}{y}}\le \norm{U^t x} \norm{y}=1,\]
where equality holds if and only if $U^t x =\gamma y$ for some $\gamma \in \cx$. Thus, an equivalent definition of perfect state transfer is that for some $t\in\ints$,
\[\abs{\ip{U^t x}{y}}=1.\]

The same generic walk is said to admit \textsl{pretty good state transfer} from $x$ to $y$, if there exists a unimodular $\gamma \in \cx$ such that, for any $\epsilon >0$, there is a time $t\in\ints$ with 
\[\norm{U^tx - \gamma y}<\epsilon.\]
Similar to perfect state transfer, pretty good state transfer can be defined using $\abs{\left<U^t x, y\right>}$.

\begin{lemma}\label{lem:pgstdef2}
Pretty good state transfer occurs from $x$ to $y$ if and only if for any $\epsilon>0$, there is $t\in\ints$ such that
\[\abs{\ip{U^t x}{ y}}>1-\epsilon.\]
\end{lemma}
\proof
Suppose first that pretty good state transfer occurs. Then there exists $\gamma\in\cx$ such that for any $\epsilon>0$, there is a time $t\in\ints$ with
\[\norm{U^tx - \gamma y}<\epsilon,\]
that is, 
\[U^t x =\gamma y +z,\]
where $z$ is a vector with norm less than $\epsilon$. Thus,
\[\abs{\ip{U^t x}{y}}  = \abs{\gamma\ip{y}{y}+\ip{y}{z}}\ge 1-\abs{\ip{y}{z}}>1-\epsilon.\]
We now prove the backward direction. Given $\delta>0$, let $t_{\delta}\in \ints$ be a time such that
\[\abs{\ip{U^{t_{\delta}} x}{ y}}>1-\delta^2/2.\]
Define 
\[\gamma_{\delta}: = \frac{\abs{\ip{y}{U^{t_{\delta}} x}}}{\ip{y}{U^{t_{\delta}} x}}.\]
Then $\abs{\gamma_{\delta}}=1$. Moreover, 
\begin{align*}
\norm{U^{t_{\delta}} x - \gamma_{\delta} y}^2 &= \ip{U^{t_{\delta}}x - \gamma_{\delta}y}{U^{t_{\delta}} x - \gamma_{\delta} y } \\
&= 2- 2\mathrm{Re}(\gamma_{\delta} \ip{y} {U^{t_{\delta}} x})\\
&= 2-2\abs{\ip{y}{U^{t_{\delta}} x}}\\
&<\delta^2.
\end{align*}
Hence, there is a sequence of integer times $\{t_k\}$ and a sequence of unimodular complex numbers $\{\gamma_k\}$ such that
\[\lim_{k\to\infty} \norm{U^{t_k} x - \gamma_k y}=0.\]
Since the unit circle is compact, $\{\gamma_k\}$ has a limit point $\gamma$ with $\abs{\gamma}=1$. Therefore, for any $\epsilon>0$, there exists $k$ such that 
\[\norm{U^{t_k} x - \gamma_k y}<\frac{\epsilon}{2},\quad \norm{\gamma_k-\gamma}<\frac{\epsilon}{2},\]
and so
\[\norm{U^{t_k} x- \gamma y }\le \norm{U^{t_k} x - \gamma_k y} + \norm{\gamma_k y-\gamma y}<\epsilon.\tag*{\sqr53}\]

We can say more about these two phenomena using the spectral information of $U$. Since $U$ is unitary, its eigenvalues are of the form $e^{i\theta_r}$ for some $\theta_r \in \re$. Let $F_r$ denote the orthogonal projection onto the $e^{i\theta_r}$-eigenspace. Then $U$ has spectral decomposition
\[U=\sum_r e^{i\theta_r} F_r.\]
This allows us to derive three inequalities for $\abs{\ip{U^t x}{y}}$. 

\begin{lemma}\label{lem:3ineq}
Let $U$ be a unitary matrix with spectral decomposition 
\[U=\sum_r e^{i\theta_r} F_r.\]
For any time $t\in\ints$ and any two states $x$ and $y$, we have
\begin{align}
\abs{\ip{U^t x}{ y}}&\le \sum_r \abs{ \ip{F_r x} {F_r y}}\label{eqn:triangle}\\
&\le \sum_r \norm{F_r x} \norm {F_r y} \label{eqn:CS1}\\
&\le\sqrt{\sum_r \norm{F_r x}^2} \sqrt{\sum_r \norm{F_r y}^2}\label{eqn:CS2}\\
&= 1 \label{eqn:1}
\end{align}
Moreover, \eqref{eqn:CS1} is tight if and only if $F_r x$ and $F_r y$ are parallel for each $r$, and \eqref{eqn:CS2} is tight if and only if $\norm{F_r x} = \norm{F_r y}$ for each $r$.
\end{lemma}
\proof
First, note that
\[\left<U^t x, y\right>=\left<\sum_r e^{i t \theta_r} F_r x, y\right>= \sum_r  e^{it\theta_r} \left< F_r x, F_r y\right>.\]
Now take the absolute value on both sides. Inequality \eqref{eqn:triangle} follows from the triangle inequaity, Inequalities \eqref{eqn:CS1} and \eqref{eqn:CS2} follow from Cauchy Schwarz inequality, and Equation \eqref{eqn:1} follows from the fact that $x$ and $y$ are unit vectors. 
\qed

From Lemma \ref{lem:3ineq}, we see that \eqref{eqn:CS1} and \eqref{eqn:CS2} are simultaneously tight if and only if 
\begin{equation}
F_r x = \mu_r F_r y, \quad \abs{\mu_r}=1, \quad \forall r \label{eqn:sc}
\end{equation}
Following the notion in continuous-time quantum walks, we say two states $x$ and $y$ are \textsl{strongly cospectral} if they satisfy \eqref{eqn:sc}. Clearly, strong cospctrality is a necessary condition for perfect state transfer. The next lemma shows that it is also necessary for pretty good state transfer.

\begin{lemma}
If pretty good state transfer occurs from $x$ to $y$, then $x$ and $y$ are strongly cospectral. 
\end{lemma}
\proof
Suppose $x$ and $y$ are not strongly cospectral. Then at least one of \eqref{eqn:CS1} and \eqref{eqn:CS2} in Lemma \ref{lem:3ineq} is a strict inequality. Thus, 
\[\sum_r \abs{\ip{F_r x}{F_ry}}<1,\]
and there exists $\epsilon>0$ such that
\[\epsilon<1-\sum_r \abs{\ip{F_r x}{F_ry}} .\]
Now for any $t\in\ints$,
\[\abs{\ip{U^t x} {y}} <1-\epsilon.\]
It follows from Lemma \ref{lem:pgstdef2} that pretty good state transfer does not occur from $x$ to $y$.
\qed

We give a spectral characterization for perfect and pretty good state transfer. To start, we define
\[\Theta_x:=\{r: F_r x \ne 0\}\]
and call it the \textsl{eigenvalue support} of the state $x$. It is clear that strongly cospectral states have the same eigenvalue support.

\begin{theorem}\label{thm:pst}
Let $U$ be a unitary matrix with spectral decomposition
\[U = \sum_r e^{i\theta_r} F_r.\]
Perfect state transfer occurs from $x$ to $y$ if and only if both of the following hold.
\begin{enumerate}[(i)]
\item $x$ and $y$ are strongly cospectral, that is, for each $r$, there is a unimodular $\mu_r\in\cx$ such that 
\[F_r x = \mu_r F_r y.\]
\item There is $t\in\ints$ such that for all $r, s\in\Theta_x$,
\[e^{it\theta_r} \mu_r = e^{it\theta_s} \mu_s.\]
\end{enumerate}
\end{theorem}
\proof
Note that
\[U^t x = \gamma y\]
if and only if for each $r$,
\[e^{it\theta_r} F_r x= \gamma F_r y.\]
Thus, perfect state transfer implies (i) and (ii). Conversely, suppose (i) and (ii) holds. Then we may let $\gamma=e^{it\theta_r} \mu_r$ for some $r\in \Theta_x$. It follows that $U^t x= \gamma y$.
\qed

\begin{theorem}\label{thm:pgst}
Let $U$ be a unitary matrix with spectral decomposition
\[U = \sum_r e^{i\theta_r} F_r.\]
Pretty good state transfer occurs from $x$ to $y$ if and only if both of the following hold.
\begin{enumerate}[(i)]
\item $x$ and $y$ are strongly cospectral, that is, for each $r$, there is a unimodular $\mu_r\in\cx$ such that 
\[F_r x = \mu_r F_r y.\]
\item There is a unimodular $\gamma\in\cx$ such that, for any $\epsilon>0$, there exists $t\in\ints$ such that for each $r\in \Theta_x$, 
\[\abs{e^{it\theta_r} \mu_r- \gamma }<\epsilon.\]
\end{enumerate}
\end{theorem}
\proof
It suffices to show that given (i), pretty good state transfer is equivalent to (ii). Indeed,
\[U^t x -\gamma y = \sum_{r\in \Theta_x} (e^{it\theta_r} \mu_r-\gamma) F_r y,\]
and so 
\begin{align*}
\norm{U^t x -\gamma y}^2 &= \ip{\sum_{r\in \Theta_x} (e^{it\theta_r} \mu_r-\gamma) F_r y}{\sum_{r\in \Theta_x} (e^{it\theta_r} \mu_r-\gamma) F_r y}\\
&=\sum_{r\in\Theta_x} \abs{e^{it\theta_r} \mu_r - \gamma}^2 \norm{F_ry}^2.
\end{align*}
Therefore, (ii) holds if and only if for any $\epsilon>0$,  there is $t\in\ints$ such that
\[\norm{U^t x -\gamma y}<\epsilon.\tag*{\sqr53}\]

\section{A special form of strong cospectrality}
In this section, we study pretty good state transfer with a spectral form of strong cospectrality, that is,
\[F_r x = \mu_r F_r y, \quad \forall r \in \Theta_x,\]
where $\mu_r$ is some root of unity. For shortness, we say $x$ and $y$ are \textsl{$m$-strongly cospectral} if for some positive integer $m$,
\[F_r x = e^{2\pi i \sigma_r/m} F_r y, \quad \forall r\in\Theta_x,\]
where $\sigma_r \in \ints_m$ and at least one of them is coprime to $m$. The following version of Kronecker's Theorem turns out to be useful.


\begin{theorem}[Kronecker]\label{thm:Kronecker}
Given $\seq{\alpha}{1}{2}{n}\in \re$ and $\seq{\beta}{1}{2}{n}\in \re$, the following are equivalent.
\begin{enumerate}[(i)]
\item For any $\epsilon>0$, the system
\[\abs{q\alpha_r - \beta_r-p_r}<\epsilon, \quad r=1,2,\cdots, n\]
has a solution $\{q, \seq{p}{1}{2}{n}\} \in \ints^{n+1}$.
\item For any set $\{\seq{\ell}{1}{2}{n}\}$ of integers such that 
\[\ell_1 \alpha_1 + \cdots + \ell_n \alpha_n \in \ints,\]
we have
\[\ell_1 \beta_1 + \cdots + \ell_n \beta_n \in \ints.\]
\end{enumerate}
\end{theorem}

We now give a number theoretic characterization for pretty good state transfer with $m$-strongly cospectral states.

\begin{theorem}\label{thm:abstract_pgst}
Let $U$ be a unitary matrix with spectral decomposition
\[U = \sum_r e^{i\theta_r} F_r.\]
Let $x$ and $y$ be two states that are $m$-strongly cospectral. Pretty good state transfer occurs from $x$ to $y$ if and only if for any set of integers $\{\ell_r: r\in\Theta_x\}$ such that 
\[\sum_{r\in\Theta_x}\ell_r \theta_r\equiv 0\pmod{2\pi},\quad \sum_{r\in\Theta_x}\ell_r =0,\]
we have
\[\sum_{r\in\Theta_x} \ell_r \sigma_r \equiv 0\pmod{m}\]
\end{theorem}

\proof
By Theorem \ref{thm:pgst}, pretty good state transfer occurs if and only if there exists $\gamma\in\cx$ such that, for any $\epsilon>0$, there is a time $t\in \ints$ with
\[\abs{e^{i(t\theta_r+2\pi \sigma_r/m)}-\gamma}<\epsilon,\quad \forall r;\]
this is equivalent to the existence of $\delta\in\re$ such that for any $\epsilon>0$, there is $t\in\ints$ with
\[\abs{t \frac{\theta_r}{2\pi} - \left(\frac{\delta}{2\pi} - \frac{\sigma_r}{m}\right) - p_r}<\epsilon.\]
Now apply Theorem \ref{thm:Kronecker} with
\[\alpha_r = \frac{\theta_r}{2\pi},\quad \beta_r = \frac{\delta}{2\pi} - \frac{\sigma_r}{m}.\]
We see that pretty good state transfer occurs if and only if the following holds:
\begin{enumerate}[(a)]
\item There exists $\delta\in \re$ such that, for any set of integers $\{\ell_r: r\in\Theta_x\}$ such that 
\[\sum_{r\in\Theta_x}\ell_r \theta_r \equiv 0 \pmod{2\pi},\]
it holds that
\begin{equation}
\delta \sum_{r\in\Theta_x} \ell_r - \frac{2\pi}{m} \sum_{r\in\Theta_x} \ell_r \sigma_r \equiv 0\pmod{2\pi}.\label{eqn:delta}
\end{equation}
\end{enumerate}
We claim that (a) is equivalent to (b):
\begin{enumerate}[(b)]
\item For any set of integers $\{\ell_r\}$ such that 
\begin{equation}
\sum_{r\in\Theta_x} \ell_r \theta_r \equiv 0 \pmod{2\pi},\label{eqn:assump}
\end{equation}
if $\sum_{r\in\Theta_x} \ell_r =0$, then 
\begin{equation}
\sum_{r\in\Theta_x} \ell_r \sigma_r\equiv 0\pmod{m}.\label{eqn:conc}
\end{equation}
\end{enumerate}
The direction (a) $\implies$ (b) is obvious. For the direction (b) $\implies$ (a), consider all sets of integers $\{\ell_r:r\in\Theta_x\}$ satisfying \eqref{eqn:assump}. There are two possibilities. 
\begin{enumerate}[(i)]
\item All such sets satisfy \eqref{eqn:conc}. Then $\delta=0$ is a common solution to \eqref{eqn:delta}.
\item Some set does not satisfy \eqref{eqn:conc}. Then
\[\left\{\sum_{r\in\Theta_x} \ell_r: \sum_{r\in\Theta_x} \ell_r \theta_r \equiv 0\pmod{2\pi}\right\}\]
forms a non-trivial additive subgroup of $\ints$. Thus, it is cyclic and has a generator $g$. Pick $\{\ell_r\}$ satisfying \eqref{eqn:delta} and
\[\sum_{r\in\Theta_x} \ell_r=g.\]
We claim that 
\[\delta = \frac{2\pi \sum_{r\in\Theta_x} \ell_r \sigma_r }{mg},\]
is a common solution to \eqref{eqn:delta}. To see this, pick any set of integers $\{\ell_r'\}$ satisfying \eqref{eqn:assump}. Define a new set of integers by
\[\eta_r:=  \frac{\sum_{s\in\Theta_x} \ell_s'}{g} \ell_r-\ell_r'.\]
Then 
\[\quad \sum_{r\in\Theta_x} \eta_r \theta_r \equiv 0\pmod{2\pi}, \quad \sum_{r\in\Theta_x} \eta_r =0.\]
Thus by (b), 
\[\sum_{r\in\Theta_x} \eta_r \sigma_r \equiv 0\pmod{m}.\]
Therefore,
\begin{align*}
\delta \sum_{r\in\Theta_x} \ell_r' - \frac{2\pi}{m} \sum_{r\in\Theta_x} \ell_r'\sigma_r&\equiv \frac{2\pi \sum_{s:s\in\Theta_x}\ell_s\sigma_s}{mg}\sum_{r\in\Theta_x} \ell_r' - \frac{2\pi}{m} \sum_{r\in\Theta_x} \ell_r'\sigma_r\\
&\equiv \frac{2\pi}{m} (\sum_{r\in\Theta_x} \eta_r \sigma_r)\\
&\equiv 0\pmod{m} \tag*{\sqr53}
\end{align*}
\end{enumerate}

Note that the condition in Theorem \ref{thm:abstract_pgst} is symmetric about $x$ and $y$. 

\begin{corollary}
Let $x$ and $y$ be $m$-strongly cospectral states. If there is pretty good state transfer from $x$ to $y$, then there is pretty good state transfer from $y$ to $x$.
\end{corollary}

\section{Quantum walk with respect to $(X, W)$}
We apply Theorem \ref{thm:abstract_pgst} to the quantum walk defined in Section \ref{sec:dqw}. Let $X$ be a graph and $W$ a normalized complex adjacency matrix. Recall that 
\[U = R(2N_t^*N_t-I),\]
where $R$ is the arc-reversal permutation matrix, and $N_t$ is the normalized tail-arc incidence with $(N_t)_{a,(a,b)} = w_{ab}$. It was shown in Section \ref{sec:spd} that the spectral decomposition of $U$ is largely determined by the following Hermitian adjacency matrix of $X$:
\[H = W \circ W^*.\]
By Lemma \ref{lem:paral}, strongly cospectrality of states $N_t^*e_a$ and $N_t^* e_b$, which is a property of the eigenprojections of $U$, translates into a property of the eigenprojections of $H$.

For convenience, we will extend a few definitions for states to vertices. Let the spectral decomposition of $H$ be 
\[H = \sum_{\lambda} \lambda E_{\lambda}.\]
The \textsl{eigenvalue support} of a vertex $a$, denoted $\Lambda_a$, is the set
\[\Lambda_a:=\{\lambda: E_{\lambda} e_a \ne 0\}.\]
We say two vertices $a$ and $b$ are \textsl{strongly cospectral} if for each eigenvalue $\lambda\in\Lambda_a$, there is a unimodular $\alpha_{\lambda}\in\cx$ such that 
\[E_{\lambda} e_a = \alpha_{\lambda} E_{\lambda} e_b;\]
if, in addition, each $\alpha_r$ is an $m$-th root of unity and at least one of them is primitive, then $a$ and $b$ are \textsl{$m$-strongly cospectral}. The following result is a direct consequence of Lemma \ref{lem:paral}.

\begin{lemma}
Two vertices $a$ and $b$ are ($m$-)strongly cospectral if and only if the states $N_t^*e_a$ and $N_t^*e_b$ are  ($m$-)strongly cospectral.
\end{lemma}

We devote the rest of the paper to perfect and pretty good state transfer between $m$-strongly cospectral vertices. Let 
\[\Lambda_{ab}^k:=\{\lambda \in \Lambda_a: E_{\lambda} e_a = e^{2\pi i k/m} E_{\lambda} e_b\}.\]
It is clear that
\[\Lambda_a = \Lambda_b =\sqcup_k\Lambda_{ab}^k\]
for $m$-strongly cospectral vertices $a$ and $b$. Moreover, there are at least two classes in this partition, as otherwise we would have
\[e_a = \sum_{\lambda\in\Lambda_a} E_{\lambda} e_a = e^{2\pi i k/m} \sum_{\lambda\in\Lambda_a} E_{\lambda} e_b=e^{2\pi i k/m} e_b.\] 
Our next characterizations is based on this partition of eigenvalue support.

\begin{theorem}\label{thm:char}
Let $a$ and $b$ be two vertices that are $m$-strongly cospectral. Let $\beta_{\pm 1}=1$ if $\pm 1\in \Lambda_a$, and $\beta_{\pm 1}=0$ otherwise. There is pretty good state transfer between $N_t^*e_a$ and $N_t^*e_b$ if and only if for any set of integers
\[\{\ell_{\lambda}: \lambda\in \Lambda_a\backslash\{1,-1\}\} \cup\{\ell_{\lambda}': \lambda\in \Lambda_a\backslash\{1,-1\}\}\cup\{\ell_1, \ell_{-1}\},\]
the two conditions
\begin{equation}\label{eqn:pgst1}
\sum_{\lambda\in \Lambda\backslash \{\pm1\}} (\ell_{\lambda}-\ell_{\lambda}') \arccos(\lambda) + \ell_{-1} \beta_{-1} \pi \equiv 0 \pmod{2\pi}
\end{equation}
and
\begin{equation}\label{eqn:pgst2}
\sum_{\lambda\in \Lambda\backslash \{\pm1\}} (\ell_{\lambda}+\ell_{\lambda'}) + \ell_1\beta_1 + \ell_{-1}\beta_{-1} =0
\end{equation}
imply
\begin{equation}\label{eqn:pgst3}
\sum_{k\in\ints_m} k\left(\sum_{\lambda\in \Lambda_{ab}^k} (\ell_{\lambda} + \ell_{\lambda}') \right)+ \ell_1\beta_1\sigma_1 + \ell_{-1} \beta_{-1} \sigma_{-1}\equiv 0\pmod{m}.
\end{equation}
\end{theorem}
\proof
By Theorem \ref{thm:eprojs}, each eigenvalue $\lambda$ of $H$ contributes to one or two eigenvalues $e^{\pm i\theta} $ of $U$, where $\theta = \arccos(\lambda)$. The statement now follows from Theorem \ref{thm:abstract_pgst}.
\qed

Theorem \ref{thm:char} imposes a strong condition on the structure of the eigenvalue support.

\begin{theorem}\label{thm:m}
Let $a$ and $b$ be two vertices that are $m$-strongly cospectral. If pretty good state transfer occurs between $N_t^*e_a$ and $N_t^*e_b$, then $m$ is even, and there is $k\in\ints_m$ such that
\[\Lambda_a = \Lambda_{ab}^k \sqcup \Lambda_{ab}^{m/2+k}.\]
\end{theorem}
\proof
Pick $k\in\ints_k^*$ for which $\Lambda_{ab}^k$ is non-empty, and let $\mu\in\Lambda_{ab}^k$. Suppose there exists 
\begin{equation}
n\in\ints_m\backslash\{k, m/2+k\}\label{eqn:nassump}
\end{equation}
for which $\Lambda_{ab}^n$ is also non-empty. Let $\eta\in\Lambda_{ab}^n$. We define a set of integers 
\[\{\ell_{\lambda}: \lambda\in \Lambda_a\backslash\{1,-1\}\} \cup\{\ell_{\lambda}': \lambda\in \Lambda_a\backslash\{1,-1\}\}\cup\{\ell_1, \ell_{-1}\}\]
based on $\eta$. If $\eta=1$ or $\eta=-1$, we choose
\[\ell_{\mu} = \ell_{\mu}'=1,\quad \ell_{\eta}=-2, \quad \mathrm{and} \quad \ell_{\lambda} = \ell_{\lambda}'=0, \quad \forall \lambda\ne \mu,\eta.\]
and if $\eta \ne \pm 1$, we choose
\[\ell_{\mu}=\ell_{\mu}'=1,\quad \ell_{\eta}=\ell_{\eta}'=-1, \quad \mathrm{and}\quad \ell_{\lambda} = \ell_{\lambda}'=0,\quad \forall \lambda\ne \mu,\eta.\]
Then Equations \eqref{eqn:pgst1} and \eqref{eqn:pgst2} hold. Now, the left hand side of Equation \eqref{eqn:pgst3} becomes $2(k-n)$, which is divisible by $m$ only when $m=2$ or $k-n=m/2$, both contradicting the assumption \eqref{eqn:nassump}.
\qed
 
 We now give a characterization for perfect state transfer, which is a special case of pretty good state transfer, for $m$-strongly cospectral vertices. 

\begin{theorem}\label{thm:pstchar}
Let $m$ be an even positive integer, and let $a$ and $b$ be two vertices that are $m$-strongly cospectral. Perfect state transfer occurs at time $t$ between $N_t^*e_a$ and $N_t^*e_b$ if and only if the following hold.
\begin{enumerate}[(i)]
\item There is $k\in\ints_m$ such that
\[\Lambda_a = \Lambda_{ab}^k \sqcup \Lambda_{ab}^{m/2+k}.\]
\item For each eigenvalue $\lambda\in(-1,1)$ of $H$, 
\[t\arccos(\lambda)\]
is a multiple of $\pi$.
\item For each pair of eigenvalues $\lambda, \mu\in \Lambda_{ab}^k$ or $\lambda, \mu \in \Lambda_{ab}^{m/2+k}$, 
\[t(\arccos(\lambda) - \arccos(\mu))\]
is an even multiple of $\pi$.
\item For each pair of eigenvalues $\lambda \in \Lambda_{ab}^k$ and $\mu \in \Lambda_{ab}^{m/2+k}$, 
\[t(\arccos(\lambda) - \arccos(\mu))\]
is an odd multiple of $\pi$.
\end{enumerate}
\end{theorem}
\proof
Condition (i) in Theorem \ref{thm:pst} is equivalent to the existence of $k\in\ints_m$ such that
\[\Lambda_a = \Lambda_{ab}^k \sqcup \Lambda_{ab}^{m/2+k}.\]
For Condition (ii) in Theorem \ref{thm:pst}, note that each eigenvalue $\lambda$ of $H$ contributes to one of two eigenvalues $e^{\pm i \theta}$ where $\cos(\theta) = \lambda$. Hence, for $\lambda_r, \lambda_s \in \Lambda_{ab}$,  the condition is equivalent to
\begin{align*}
e^{it\arccos(\lambda_r)} &= e^{it\arccos(\lambda_s)},\quad \lambda_r, \lambda_s \in \Lambda_{ab}^k \text{ or } \lambda_r, \lambda_s \in \Lambda_{ab}^{m/2+k} \\
e^{it\arccos(\lambda_r)} &= -e^{it\arccos(\lambda_s)},\quad \lambda_r\in \Lambda_{ab}^k, \quad \lambda_s \in \Lambda_{ab}^{m/2+k} \\
e^{it\arccos(\lambda_r)} & = e^{-it\arccos(\lambda_r)},\quad \lambda_r\in(-1,1). \tag*{\sqr53}
\end{align*}

\section{Constructions}
We have reduced a problem on the quantum walk to a problem on certain Hermitian adjacency matrix $H$ of the underlying graph. However, to recover a walk from $H$, we need to find $W$ so that $H = W\circ W^*$ and $W\circ \comp{W}$ is row-stochastic. 
\begin{lemma}\label{lem:recover}
Let $H=(h_{ab})_{ab}$ be a Hermitian adjacency matrix of a graph such that for every vertex $a$,
\[\sum_{b\sim a} \abs{h_{ab}}=1.\]
For each edge $\{a,b\}$, pick a direction $(a,b)$, a unimodular complex number $\delta_{ab}$, and let
\[w_{ab} = \sqrt{\abs{h_{ab}}} \delta_{ab},\quad w_{ba} = \frac{\comp{h_{ab}}}{\comp{w_{ab}}}.\]
Let $W=(w_{ab})_{ab}$. Then $H = W \circ W^*$ and $W\circ \comp{W}$ is row-stochastic.
\end{lemma}
\proof
One can verity that
\[w_{ab} \comp{w_{ba}} = h_{ab}\]
and
\[\sum_{b\sim a} w_{ab} \comp{w_{ab}}=\sum_{b\sim a} \abs{h_{ab}}=1.\tag*{\sqr53}\]

Lemma \ref{lem:recover} finds more than one choice of $W$ for the same $H$. From a physical point of view, this means we can engineer pretty good state transfer of different states on the same graph, by assigning different phases $\delta_{ab}$.

A real number $\theta$ is called a \textsl{pure geodetic angle} if any one of its six squared trigonometric functions is either rational or infinite. From Theorem \ref{thm:char} Equation \eqref{eqn:pgst1}, we see that pretty good state transfer is related to linear dependence of certain angles over the rationals. Such linear relations have been well-studied for geodetic angles; we cite two useful results.

\begin{theorem}\cite[Ch 11]{Bergen2009}\label{thm:tan_rat_pi} %
If $q\in \rats$ and $q\pi$ is a pure geodetic angle, then 
\[\tan(q\pi) \in \left\{0, \pm{\sqrt{3}}, \pm \frac{1}{\sqrt{3}}, \pm 1\right\}.\]
\end{theorem}

\begin{theorem}\cite{Conway1999}\label{thm:rat_lb_geo}
If the value of a rational linear combination of pure geodetic angles is a rational multiple of $\pi$, then so is the value of its restriction to those angles whose tangents are rational multiple of any given square root.
\end{theorem}

We construct several infinite families of quantum walks that admit pretty good state transfer.

\subsection{Normalized adjacency matrix}
The first two families are with respect to $(X,W)$, where $X$ is distance regular and $W\circ W^*$ is the normalized adjacency matrix of $X$. By Lemma \ref{lem:recover}, more than one choices of $W$ will work; the simples example is $W=(w_{ab})_{ab}$ where $w_{ab}=1/\sqrt{\deg(a)\deg(b)}$.

\begin{theorem}
Let $X$ be the cocktail party graph $\comp{n K_2}$, and let $H$ be the normalized adjacency matrix of $X$. Pick any $W$ such that $H=W\circ W^*$ and $W\circ \comp{W}$ is row-stochastic. The quantum walk with respect to $(X, W)$ has pretty good state transfer between the antipodal vertices.
\end{theorem}
\proof 
The cases where $n=2$ and $n=3$ were proved in \cite{Zhan2019}; in fact, $\comp{2K_2}$ and $\comp{3K_2}$ have perfect state transfer. We show that pretty good state transfer occurs for larger $n$. Let $H$ be the normalized adjacency matrix of $X$. Since $X$ is regular with valency $2n-2$, 
\[H = \frac{1}{2n-2} A(X),\]
and so strong cospectrality relative to $H$ is equivalent to strong cospectrality relative to $A(X)$. It is known (see, for example, \cite{Coutinho2015a}) that the antipodal vertices $a$ and $b$ of $X$ are $2$-strongly cospectral, and, relative to $A(X)$,
\[\Lambda_{ab}^0=\{2n-2, -2\},\quad \Lambda_{ab}^1 = \{0\}.\]
Now suppose
\[(\ell_{-2}-\ell_{-2}') \arccos\left(-\frac{1}{n-1}\right)+(\ell_0-\ell_0') \frac{\pi}{2} \equiv 0\pmod{2\pi}.\]
By Theorem \ref{thm:tan_rat_pi}, $\arccos(-1/(n-1))$ is linearly independent with $\pi$ over the rationals. Hence $\ell_{-2}=\ell_{-2}'$, and so $\ell_0-\ell_0'$ is even. It follows from Theorem \ref{thm:char} and
\[\ell_0+\ell_0'\equiv 0 \pmod{2}\]
that pretty good state transfer occurs on $X$.
\qed

To construct the second family, we derive some tools to treat rational eigenvalue support.

\begin{lemma}\label{lem:sqff}
Let $d,\seq{\lambda}{1}{2}{k}$ be positive integers. Suppose 
\begin{enumerate}[(i)]
\item $\lambda_r\in(0, d/2) \cup (d/2, d)$ for each $r$, and
\item the square free parts of
\[d^2-\lambda_1^2,\quad d^2-\lambda_2^2, \quad \cdots,\quad d^2-\lambda_k^2\]
are pairwise distinct.
\end{enumerate}
Then the angles 
\[\pi,\quad \arccos(\lambda_1/d),\quad \arccos(\lambda_2/d),\quad \cdots, \quad \arccos(\lambda_k/d)\]
are linearly independent over $\rats$.
\end{lemma}
\proof
Let $\theta_r = \arccos(\lambda_r/d)$. Then $\tan(\theta_r)$ is a rational multiple of $\sqrt{d^2-\lambda_r^2}$. Suppose for some $q_r\in \rats$,
\[q_1 \theta_1 + q_1 \theta_2 + \cdots+ q_k \theta_k\in \rats\pi.\]
By Theorem \ref{thm:rat_lb_geo}, each $\theta_r$ is a rational multiple of $\pi$. It follows from Theorem \ref{thm:tan_rat_pi} that $\tan(\theta_r)\in\{0, \pm \sqrt{3}, \pm 1/\sqrt{3}, \pm1\}$. However, this is impossible as $\lambda_r$ is an integer and $\lambda_r \ne d/2$.
\qed

\begin{lemma}\label{lem:prime}
Let $p$ be an odd prime. For any two distinct integers $j,k\in\{1,2,\cdots, (p-1)/2\}$, the square free parts of
\[j(p-j),\quad k(p-k)\]
are distinct.
\end{lemma}
\proof 
First note that two integers $m$ and $n$ have the same square free part if and only if $mn$ is a perfect square. Now write
\[p=a+b = c+d\]
for some positive integers $a, b, c, d$. Suppose $ac$ and $bd$ have the same square free part. Then $abcd$ is a perfect square. It follows that $ac$ and $bd$ have the same square free part, say $\Delta_1$, and $ad$ and $bc$ have the same square free part, say $\Delta_2$. As
\[p^2 = (a+b)(c+d) = (ac+bd)+(ad+bc),\]
$\gcd(\Delta_1, \Delta_2)$ divides $p^2$. If $p$ were a factor of $\Delta_1$, it would divide $a$ or $c$, which contradicts our choice of $a, b, c, d$. Hence, $\Delta_1$ and $\Delta_2$ are coprime, and so $\gcd(ac, ad)$ is a perfect square. On the other hand, $c$ and $d$ are also coprime, which implies
\[\gcd(ac, ad) = a \gcd(c,d) = a.\]
Therefore, $a$ is a perfect square. A similar argument shows that $b$, $c$ and $d$ are all perfect squares. Since every prime can be expressed in at most one way as a sum of two squares, either $a=c$ or $a=d$.
\qed

We now show that pretty good state transfer occurs on every hypercube $Q_p$ with prime $p$.

\begin{theorem}
Let $p$ be a  prime. Let $X$ be the hypercube $Q_p$, and let $H$ be the normalized adjacency matrix of $X$. Pick any $W$ such that $H=W\circ W^*$, and $W\circ \comp{W}$ is row-stochastic. The quantum walk with respect to $(X, W)$ has pretty good state transfer between the antipodal vertices.
\end{theorem}
\proof
We prove this for odd prime $p$, as $Q_2 = \comp{2K_2}$. It is known \cite{Coutinho2015a} that the adjacency eigenvalues of $Q_p$ are
\[\lambda_r = p - 2r, \quad r= 0,1,\cdots,p,\]
and, relative to $A(X)$, the antipodal vertices $a$ and $b$ are $2$-strongly cospectral with
\[\Lambda_{ab}^1 = \{p - 2r: r\text{ odd}\},\quad \Lambda_{ab}^0 = \{p - 2r: r\text{ even}\} \]
Let $\theta_r =\arccos(\lambda_r/p)$. Then $\tan(\theta_r)$ is a rational multiple of $\sqrt{(p-r) r}$. By Lemma \ref{lem:sqff}, Lemma \ref{lem:prime} and symmetry of the eigenvalues,
\[\pi, \theta_1, \theta_3, \cdots,\theta_p\]
are linearly independent over $\rats$. Now 
\[\sum_{0<r<p} (\ell_r - \ell_r') \theta_r  + \ell_{-p} \pi\equiv 0 \pmod{2\pi}.\]
We may rewrite the left hand side as
\[\sum_{0<r<p,  \text{ $r$ odd}}(\ell_r - \ell_r' -\ell_{p-r}+\ell_{p-r}') \theta_r  + \left(\ell_{-p} + \sum_{0<r<p,  \text{ $r$ odd}}(\ell_{p-r} -\ell_{p-r}')\right) \pi.\]
It follows that
\[\ell_{p-r}-\ell_{p-r}'=\ell_r-\ell_r'\]
for every odd $r$, and 
\[\ell_{-p}+\sum_{0<r<p,  \text{ $r$ odd}}\left(\ell_r + \ell_r'\right)\equiv 0 \pmod{2},\]
 Thus, Equation \eqref{eqn:pgst3} in Theorem \ref{thm:char} holds.
\qed

It remains open to determine $n$ for which the quantum walk relative to $(Q_n, A(Q_n)/\sqrt{n})$ has pretty good state transfer.


\subsection{Cyclic covers}
A graph $Y$ is called a \textsl{cover} of $X$ if there is a surjective homomorphism from $Y$ to $X$ that is locally bijective. For instance, the hypercube $Q_3$ covers the complete graph $K_4$, where the covering map sends each pair of antipodal vertices in $Q_3$ to a distinct vertex in $K_4$.

In general, covers can be constructed through \textsl{symmetric arc functions}, that is, a function $f$ from the arcs of $X$ to $\sym{r}$, such that 
\[f(a,b) = f(b,a)^{-1}\]
for any two adjacent vertices $a$ and $b$. Let $Y^f$ be the graph with vertex set $V(X) \times \{1,2,\cdots,r\}$ and edge set
\[\{((a, i), (b, f(a, b)f(i))): a\sim b\}\]
Then $Y^f$ is a cover of $X$, called an \textsl{$r$-fold cover}. An $r$-fold cover is \textsl{cyclic} if the image of $f$ generates a cyclic group.

If $H$ is a Hermitian adjacency matrix of $X$ whose entries are $r$-th roots of unity, then $H$ naturally defines a symmetric arc function and hence an $r$-fold cover of $X$. Moreover, the spectrum of $H$ is a subset of the adjacency spectrum of the cover. An important topic in algebraic graph theory is to construct covers of graphs whose adjacency spectrum differs very little from that of the underlying graph; we refer interested readers to \cite{Godsil1992} and \cite{Godsil2021} for more details.

Given any graph $X$, the \textsl{Seidel adjacency matrix} of $X$ is $J-I-2A$, which can been seen as a signed adjacency matrix, and hence a double cover, of the complete graph. While complete graphs on more than two vertices do not have strongly cospectral vertices relative to the adjacency matrix, they do relative to some signed adjacency matrices. We use these matrices to construct walks that admit pretty good state transfer on complete graphs. First, we cite an old result due to Erdos.

\begin{theorem}\cite{Erdos1939}\label{thm:ps}
A product of consecutive positive integers is never a perfect square.
\end{theorem}

\begin{theorem}
Let $n$ be a positive integer that is not twice a perfect square. Let $X=n K_2$ and let $H$ be the Seidel adjacency matrix of $X$. Let $W$ be any matrix such that $H = (2n-1) (W\circ W^*)$ and $W\circ \comp{W}$ is a row-stochastic matrix. The quantum walk with respect to $(K_{2n}, W)$ has pretty good state transfer between every pair of adjacent vertices in $X$.
\end{theorem}
\proof
The Seidel adjacency matrix of $X$ is
\[H = I_n\otimes I_2 + (J_n- I_n)\otimes J_2 - I_n\otimes J_2,\]
and it has spectral decomposition
\[H = (2n-3)\cdot \frac{1}{2n} (J_n \otimes J_2) + 1\cdot \left( I_n\otimes \left(I-\frac{1}{2}J_2\right) \right)- 3\cdot \left(\left(I_n-\frac{1}{n}J_n\right) \otimes \frac{1}{2}J_2\right).\]
Let $0$ and $1$ be two adjacent vertices in $X$. We see that $0$ and $1$ are $2$-stongly cospectral, with
\[\Lambda_{01}^0 = \{2n-3, -3\},\quad \Lambda_{01}^1 = \{1\}.\]
Since each row of $H$ has $2n-1$ non-zero entries, we need to apply Theorem \ref{thm:char} to 
\[\theta_1 = \arccos\left(\frac{2n-3}{2n-1}\right),\quad \theta_2 = \arccos\left(-\frac{3}{2n-1}\right),\quad \theta_3 = \arccos\left(\frac{1}{2n-1}\right).\]
Note that the tangents of these angles are, respectively, rational multiples of square roots of
\[a_1 = 2(n-1),\quad a_2 = (n+1)(n-2),\quad a_3 =n(n-1).\]
We claim that $a_1$ and $a_3$ do not have the same square free factor, and neither do $a_2$ and $a_3$. Indeed, by our choice of $n$, the product $a_1a_3$ is not a perfect square, and by Theorem \ref{thm:ps}, $a_2a_3$ is not a perfect square. Therefore, due to Theorem \ref{thm:rat_lb_geo},
\[(\ell_1-\ell_1')\theta_1  + (\ell_2-\ell_2')\theta_2 + (\ell_3-\ell_3')\theta_3 \equiv 0\pmod{2\pi}\]
implies that $(\ell_3-\ell_3')\theta_3 $ is a rational multiple of $\pi$. Moreover, $\tan(\theta_3) = 2\sqrt{n(n-1)}$, which shows that $\theta_3$ is not a rational multiple of $\pi$. Hence, $\ell_3=\ell_3'$, and so Equation \eqref{eqn:pgst3} in Theorem \ref{thm:char} holds.
\qed

We now consider a special class of cyclic $4$-fold covers, that is, Hermitian adjacency matrices with entries in $\{\pm1, \pm i\}$, associated to oriented graphs. Given an orientation of $X$, the \textsl{skew-adjacency matrix}, denoted $S(X)$, is the matrix whose $ab$-entry is $1$ if $(a,b)$ is an arc in the orientation, $-1$ if $(b,a)$ is an arc in the orientation, and $0$ otherwise. As $H=iS(X)$ is Hermitian, the eigenvalues of $S(X)$ are purely imaginary and symmetric about $0$. Moreover, if $E_{\lambda}$ is the projection onto the $\lambda$-eigenspace of $H$, then $\comp{E_{\lambda}}$ is the projection onto the $(-\lambda)$-eigenspace of $H$.

The following lemma concerns an orientation of the Cartesian product $K_2\square X$; its skew-adjacency spectrum has been studied in \cite{Cui2013}. Here, we prove this construction preserves $2$-strong cospectrality.

\begin{lemma}\label{lem:cart}
Let $S(X)$ be a skew-adjacency matrix of an oriented graph $X$. Let $H(X)=iS(Y)$. Let $Y=K_2 \square X$ and define an orientation on $Y$ by
\[H(Y) = \pmat{H(X) & iI \\ -iI & -H(X)}.\]
Then the eigenvalues of $H(Y)$ are $\pm \sqrt{\lambda^2+1}$ for each eigenvalue $\lambda$ of $H(X)$.  Moreover, any pair of $2$-strongly cospectral vertices in $X$ are also $2$-strongly cospectral in $Y$.
\end{lemma}
\proof
First note that for any eigenvalue $\lambda$ and eigenprojection $E_{\lambda}$, 
\[E_{\lambda} e_a= \alpha E_{\lambda} e_b\] 
if and only if each eigenvector $z$ for $\lambda$ satisfies $\ip{e_a}{z} = \alpha \ip{e_b}{z}$. We now consider the eigenspaces of $H(Y)$. Suppose for some $\mu\in\re$ and vectors $x,y$,
\[H(Y) \pmat{x\\y} = \mu \pmat{x\\y}.\]
Then 
\[(\mu I-H(X))x = iy,\quad (\mu I+H(X)) y = -ix,\]
from which it follows that
\[H(X)^2x=(\mu^2-1)x,\quad H(X)^2y=(\mu^2-1)y.\]
This proves the first part of the lemma. To see the second part, suppose $a$ and $b$ are $2$-strongly cospectral relative to $H(X)$. Then for each eigenvalue $\lambda$ of $H(X)$ in the eigenvalue support of $a$, either both $\pm \lambda$ lie in $\Lambda_{ab}^0$, or both $\pm \lambda$ lie in $\Lambda_{ab}^1$. Therefore, every eigenvector $\pmat{x\\ y}$ for $H(Y)$ with eigenvalue $\mu$ satisfies $\ip{e_a}{z} = \alpha \ip{e_b}{z}$, where $\alpha=1$ if $\sqrt{\mu^2-1} \in \Lambda_{ab}^0$, and $\alpha=-1$ if $\sqrt{\mu^2-1} \in \Lambda_{ab}^1$.
\qed

We show how cyclic covers may be used to generate pretty good state transfer.  Let $X$ be the line graph of the cube $Q_3$. Since there are no strongly cospectral vertices relative to $A(X)$, pretty good state transfer does not occur on the quantum walk with respect to $(X, A(X)/2)$. However, the following orientation gives a Hermitian matrix relative to which $1$ and $10$ are strongly cospectral, and together with Lemma \ref{lem:cart}, it provides an infinite family of quantum walks with pretty good state transfer.

\begin{figure}[h]
\centering
\includegraphics[width=6cm]{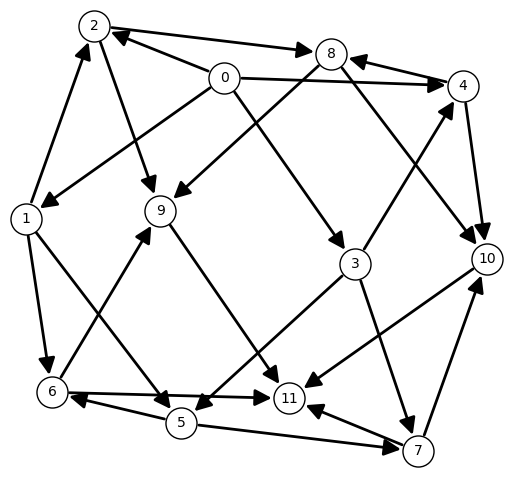}
\caption{An orientation of the line graph of $Q_3$}
\label{fig:orien}
\end{figure}

\begin{theorem}
Let $X$ be the line graph of $Q_3$ and let $S(X)$ be the skew-adjacency matrix of the orientation shown in Figure \ref{fig:orien}. Define a sequence of Hermitian matrices as follows:
\[H_4 = iS(X),\quad H_d = \pmat{H_{d-1} & iI\\ -iI & H_{d-1}}.\]
Let $X_d$ be the underlying graph of $H_d$. Pick any Hermitian matrix $W_d$ such that $H_d = d(W_d \circ W_d^*)$ and $W_d\circ \comp{W_d}$ is row-stochastic. Then for each odd prime $p$, pretty good state transfer occurs on the quantum walk with respect to $(X_p, W_p)$.
\end{theorem}
\proof
A direct computation shows that vertices $1$ and $10$ are $2$-strongly cospectral in $X$, and 
\[\Lambda_{1, 10}^0=\{\pm 2\},\quad \Lambda_{1,10}^1=\{0, \pm 12\}.\]
By Lemma \ref{lem:cart}, the same vertices are $2$-strongly cospectral in $X_d$, and, relative to $H_d$, 
\[\Lambda_{1, 10}^0=\{\pm \sqrt{d}\},\quad \Lambda_{1,10}^1=\{\pm \sqrt{d-4}, \pm\sqrt{d+8}\}.\]
Since each row of $H_d$ has exactly $d$ non-zero entries, to check pretty good state transfer we consider the following linear relation:
\[\sum_{r=1}^3 ((\ell_r-\ell_r') \theta_r  + (\ell_{-r}-\ell_{-r}')(\pi-\theta_r))\equiv 0 \pmod{2\pi},\]
where 
\[\theta_1 = \arccos\left(\frac{\sqrt{d}}{d}\right), \quad \theta_2 = \arccos\left(\frac{\sqrt{d-4}}{d}\right),\quad \theta_3 = \arccos\left(\frac{\sqrt{d+8}}{d}\right).\]
Hence the tangents of $\theta_1, \theta_2, \theta_3$ are rational multiples of square roots of the following, respectively,
\[a_1 = (d^2-d)d, \quad a_2=(d^2-(d-4))(d-4),\quad a_3=(d^2-(d+8))(d+8).\]
Now suppose $d=p$ for some odd prime $p$. By an argument similar to Lemma \ref{lem:prime}, if $a_1a_2$ or $a_1a_3$ is a perfect square, then $d$ is a perfect square, and if $a_2a_3$ is a perfect square, then both $d-4$ and $d+8$ are perfect squares. Neither of these can happen when $d=p$ is an odd prime. Thus, 
\[\ell_r - \ell_r'=\ell_{-r} - \ell_{-r}', \quad r= 1,2,3\]
and by Theorem \ref{thm:char}, pretty good state transfer occurs between vertices $1$ and $10$.
\qed

We expect more examples of pretty good state transfer arising from cyclic covers of graphs; in particular, the Hermitian adjacency matrices defined in \cite{Li2015}, \cite{Guo2015} and \cite{Mohar2020} are worth investigating.

\subsection{Real-weighted graphs}
Compared to pretty good state transfer, perfect state transfer is rarer - the eigenvalue condition in Theorem \ref{thm:pst} is particularly hard to satisfy when $H$ is the normalized adjacency matrix. In this section, we use other real weights to obtain quantum walks that admit perfect state transfer.

\begin{theorem}
Let $s$ be a positive integer and set $n=2^{s-1}-1$. There exist two sets of odd integers
\[\{\seq{p}{1}{2}{n}\}, \quad  \{\seq{q}{1}{2}{n}\}\]
such that $\gcd(p_j, q_j)=1$ for each $j$, and
\[\sum_{j=1}^n \abs{\cos\left(\frac{p_j}{q_j}\right)}<1.\]
Define
\[D = \mathrm{diag}\left(1, -1, \cos\left(\frac{p_1\pi}{q_1}\right), -\cos\left(\frac{p_1\pi}{q_1}\right), \cos\left(\frac{p_2\pi}{q_2}\right), - \cos\left(\frac{p_2\pi}{q_2}\right),\cdots\right).\]
Let $P$ be the Hadamard matrix
\[P = \pmat{1&1\\1&-1}^{\otimes s}\]
and 
\[H=PDP^T.\]
Then $H$ is a real-weighted adjacency matrix of $K_{n+1, n+1}$. Moreover, for each $W$ such that $H=(2n+2)(W\circ W^*)$ and $W\circ \comp{W}$ is row-stochastic, there is perfect  state transfer between vertices $0$ and $1$ on the quantum walk with respect to $(K_{n+1, n+1}, W)$.
\end{theorem}
\proof
Starting from $s=1$, one can construct the set 
\[\{\seq{p}{1}{2}{n}\}, \quad  \{\seq{q}{1}{2}{n}\}\]
inductively such that $p_j$ is coprime to $q_j$ and the sum 
\begin{equation}\label{eqn:sumcos}
\sum_{j=1}^n \abs{\cos\left(\frac{p_j}{q_j}\right)}<1.
\end{equation}
Now, it is easy to tell from $P$ that all vertices are $2$-strongly cospectral relative to $H$. In particular for $0$ and $1$, we have
\begin{align*}
\Lambda_{01}^0 &= \{1\} \cup \left\{\cos\left(\frac{p_j \pi} {q_j}\right): j=1,2,\cdots, n \right\},\\
\Lambda_{01}^1 &= \{-1\} \cup \left\{-\cos\left(\frac{p_j \pi} {q_j}\right): j=1,2,\cdots, n \right\}.
\end{align*}
Thus, $t=q_1q_2\cdots q_n$ satisfies the condition in Theorem \ref{thm:pst}. A bit calculation shows that $H$ has constant row sum $2n+2$ and, due to \eqref{eqn:sumcos}, all non-zero entries are positive. Therefore, there exists $W$ such that $H=(2n+2)(W\circ W^*)$ and $W\circ \comp{W}$ is row-stochastic.

\subsection{$m$-strongly cospectral vertices with $m>2$}
In all the examples we have seen so far, pretty good state transfer occurs between $2$-strongly cospectral vertices. This raises the question whether $m$ can take any value. Theorem \ref{thm:m} already rules out all odd values; we now show that pretty good state transfer can happen for any even integers $m$.

\begin{theorem}
Let $m$ be any positive even integer and let $\alpha=e^{2\pi i /m}$. Let $p>3$ be a prime. Define a Hermitian adjacency matrix $H$ of $K_4$ by
\[H = \pmat{0&  2 \comp{\alpha}& p &p\\
2\alpha & 0 & p\alpha & p\alpha\\
p & p\comp{\alpha} & 0& 2\\
p& p\comp{\alpha} & 2 & 0}.\]
For every $W$ such that $H=2(p+1)W\circ W^*$ and $W\circ \comp{W}$ is row-stochastic, there is pretty good state transfer between vertices $0$ and $1$ on the quantum walk with respect to $(K_4, W)$.
\end{theorem}
\proof
The spectral decomposition of $H$ is
\[H = PDP^*,\]
where 
\[P = \frac{1}{2}\pmat{1&1&\sqrt{2}&0\\ \alpha & \alpha & -\sqrt{2}\alpha& 0\\1 & -1&0&\sqrt{2}\\ 1 & -1 & 0 & -\sqrt{2}}\]
and
\[D =\mathrm{diag}(2+2p, 2-2p, -2, -2).\]
Hence vertices $0$ and $1$ are $m$-strongly cospectral, and 
\[\Lambda_{01}^1=\{2+2p, 2-2p\},\quad \Lambda_{01}^{m/2+1} = \{-2\}.\]
Let 
\[\theta_1 = 0,\quad \theta_2 = \arccos\left(-\frac{1-p}{1+p}\right),\quad \theta_3 = \arccos\left(-\frac{2}{1+p}\right).\]
According to Theorem \ref{thm:char}, we need to show that 
\[(\ell_2-\ell_2') \theta_2 +(\ell_3-\ell_3')\theta_3 \equiv 0\pmod{2\pi}\]
and
\[\ell_2+\ell_2'+\ell_3+\ell_3'+\ell_1=0\]
imply
\begin{equation}\label{eqn:3rd}
(\ell_2+\ell_2')+(m/2+1)(\ell_3+\ell_3')+\ell_1\equiv 0\pmod{m}.
\end{equation}
Since $\tan(\theta_1)$ is rational multiple of $\sqrt{p}$, while $\tan(\theta_3)$ is a rational multiple of $\sqrt{p^2+2p+3}$, the pure geodetic angles $\theta_1$ and $\theta_2$ are linearly independent with $\pi$ over the rationals. Thus, 
\[\ell_2 =\ell_2',\quad \ell_3=\ell_3',\]
and so \eqref{eqn:3rd} holds.
\qed

\section{Future work}
We have developed the theory for pretty good state transfer in discrete-time quantum walks. For the quantum walk with respect to $(X,W)$, the coin matrix of the quantum walk translates into a Hermitian adjacency matrix of $X$, which in turn determines the transport properties of the walk. While we have constructed some infinite families of walks that admit pretty good state transfer, this is merely a start to exploring the rich connection between the spectrum of Hermitian adjacency matrices and the behavior of quantum walks. Below we list a few questions for future research in this area.

\begin{enumerate}
\item Strong cospectrality relative to the adjacency matrix $A(X)$ has been studied in \cite{Godsil2017Smith}. Since $A(X)$ is real symmetric, any two vertices that are $m$-strongly cospectral relative to $A(X)$ must be $2$-strongly cospectral. It is worth extending the theory in \cite{Godsil2017Smith} to $m$-strong cospectrality for larger $m$, in particular relative to (i) the signed adjacency matrices, (ii) the Hermitian adjacency matrix with entries in $\{1,\pm i,0\}$ \cite{Li2015,Guo2015}, or (iii) the Hermitian adjacency matrix with entries in $\{1, e^{\pm 2\pi i/3},0\}$ \cite{Mohar2020}.

\item The basis for the rational vector space generated by geodetic angles was studied in \cite{Conway1999}. It will be interesting to apply this result to determine pretty good state transfer between vertices whose eigenvalue supports consist of integers or quadratic integers.

\item There are strongly cospectral vertices that are not $m$-strongly cospectral. For example, relative to the following Hermitian adjacency matrix \cite{Guo2015}:
\[H= \left(\begin{array}{rrrr}
0 & i & -i & i \\
-i & 0 & i & i \\
i & -i & 0 & 1 \\
-i & -i & 1 & 0
\end{array}\right),\]
 vertices $2$ and $3$ are strongly cospectral, and there is an eigenvalue $\lambda$ with 
\[E_{\lambda}e_a = e^{-i\arctan(1/2)} E_{\lambda}e_b.\]
Can pretty good state transfer occur on such type of vertices? If so, how do we test it efficiently?
\end{enumerate}

\bibliographystyle{amsplain}
\bibliography{dqw.bib}
\end{document}